\documentstyle{amsppt}
\magnification=1200
\hsize=150truemm
\vsize=224.4truemm
\hoffset=4.8truemm
\voffset=12truemm

\NoRunningHeads

\define\C{{\bold C}}

\let\thm\proclaim
\let\fthm\endproclaim

\define\pp{ P^2(\C) }
\define\p {P^3(\C) }

\newcount\tagno
\newcount\secno
\newcount\subsecno
\newcount\stno
\global\subsecno=1
\global\tagno=0
\define\ntag{\global\advance\tagno by 1\tag{\the\tagno}}

\define\sta{\ 
{\the\secno}.\the\stno
\global\advance\stno by 1}

\define\stas{\the\stno
\global\advance\stno by 1}

\define\sect{\global\advance\secno by 1
\global\subsecno=1\global\stno=1\
{\the\secno}. }

\def\nom#1{\edef#1{{\the\secno}.\the\stno}}
\def\inom#1{\edef#1{\the\stno}}
\def\eqnom#1{\edef#1{(\the\tagno)}}

\newcount\refno
\global\refno=0

\def\nextref#1{
      \global\advance\refno by 1
      \xdef#1{\the\refno}}

\def\bref {\ref\global\advance\refno by 1\key{\the\refno}}


\nextref\BER 
\nextref\BRO
 \nextref\DEM
 \nextref\SHI
 \nextref\ZAI

\topmatter

\title 
Une sextique hyperbolique dans $\p$ 
 \endtitle

\author  Julien Duval
\endauthor

\abstract\nofrills{\smc R\'esum\'e. }\ \ On construit une sextique 
hyperbolique 
dans $\p$.

\null
{\smc \hskip 2,5cm  A hyperbolic sextic surface in $\p$}

\null
\noindent
{\smc Abstract.}\ We construct a hyperbolic sextic surface in $\p$.
 \endabstract

\endtopmatter 

\document
 
 Une partie de $\p$ est {\it hyperbolique} si elle ne contient pas de 
 courbe enti\`ere, i.e. d'image holomorphe non constante 
 de $\C$. La conjecture de Kobayashi dans l'espace projectif stipule
 qu'une surface g\'en\'erique de degr\'e $\geq 5$ de $\p$ est 
 hyperbolique.
 Elle a \'et\'e d\'emontr\'ee pour les degr\'es $\geq 21$ 
 par J.-P. Demailly et J. El Goul [\DEM].
 Parall\`element, et de mani\`ere plus modeste, un certain nombre 
 d'auteurs (voir les r\'ef\'erences dans [\ZAI]) ont cherch\'e \`a construire 
 des exemples de surfaces 
 hyperboliques de degr\'e le plus bas possible. A ce jour, la 
 meilleure borne \'etait le degr\'e $8$. Le but de cette note est 
 de montrer l'existence de surfaces hyperboliques de degr\'e $6$ :
 
\thm{ Th\'eor\`eme} Il existe des sextiques hyperboliques.
\fthm

 L'argument est un raffinement de la m\'ethode de filtrage 
 (percolation) de B. 
 Shiffman et M. Zaidenberg (voir par exemple [\SHI]). 
 La question se r\'eduit \`a l'impos-sibilit\'e de faire passer une 
 courbe enti\`ere par un filtre convenable du plan projectif :
 on se ram\`ene \`a l'hyperbolicit\'e d'un compl\'ementaire de la 
 forme $\pp \setminus (Q^* \cup C^*)$, avec 
 $Q^*=Q \setminus C'$ et $C^*=C \setminus Q'$,  $Q,Q'$ \'etant
 des coniques et $C,C'$ des cubiques de $\pp$.  La souplesse de 
 choix
 de $C'$ et $Q'$ permet alors une 
 seconde r\'eduction \`a l'hyperbolicit\'e de $\pp \setminus(Q \cup 
 C)$. Celle-ci est connue si $Q \cup C$ est assez proche d'une 
 configuration plus d\'eg\'en\'er\'ee.

  D\'etaillons cette approche.

 Les deux ingr\'edients principaux en sont le lemme de Brody [\BRO] et la 
 persistance d'intersection. Le premier sera employ\'e sous la forme 
 suivante : 
 
 {\it d'une suite de courbes enti\`eres, on peut 
 extraire une sous-suite convergeant (localement uniform\'ement 
 apr\`es reparam\'etrisation) vers une courbe enti\`ere. }
 
 La seconde 
 s'\'enonce ainsi dans $\p$ : 
 
 {\it soit $(L_n)$ une suite de courbes enti\`eres 
 convergeant vers une courbe enti\`ere $L$ et $(S_n)$ une suite 
 de surfaces convergeant vers une surface $S$. Si $L$ n'est pas 
 contenue dans $S$, alors $L_n$ coupe $S_n$   
 pr\`es de tout point de $L\cap S$ pour $n$ assez grand. }
 
 \null
 
 En les combinant, on obtient, dans nombre de cas, le fait que 
 l'hyperbolicit\'e 
 est une propri\'et\'e ouverte pour les compl\'ementaires. Par exemple :
 
 \thm{Un cas d'ouverture} 
 Soient $S,S'$ deux surfaces dans $\p$ et $U$ un ouvert 
 de $S'$.
 On suppose $S\cap S'$ et $S \setminus U$ hyperboliques; 
 alors, pour toute 
 surface $S_\epsilon$ proche de $S$, $S_\epsilon \setminus U$ est 
 encore hyperbolique.  
 \fthm
 
 Sinon, on produit une suite $(L_n)$ de courbes enti\`eres 
 dans $S_{\epsilon_n} \setminus U$ convergeant vers $L$ dans $S$. Par 
 hyperbolicit\'e de $S \cap S'$, $L$ n'est pas contenue dans $S'$. 
 Par persistance d'intersection, puisque $L_n$ \'evite $U$ il en est de 
 m\^eme pour $L$. Ceci contredit l'hyperbolicit\'e de $S \setminus 
 U$. $\square$

 \null
 
 Voici maintenant le principe de filtrage de Shiffman et 
 Zaidenberg [\SHI] :
 
\thm{Filtrage} Soient $S_i=(f_i=0)$ trois surfaces de degr\'e $d_i$ 
dans $\p$ avec 
$d_3=d_1+d_2$. On 
suppose $S_1\cap S_2$, $S_1 \setminus (S_2 \setminus S_3)$
et $S_2 \setminus 
(S_1 \setminus S_3)$ hyperboliques. Alors la surface 
$S_\epsilon = (f_1f_2=\epsilon f_3)$ est hyperbolique pour 
$\epsilon$ petit non nul.  
\fthm

Sinon, on construit une suite $(L_n)$ de courbes enti\`eres 
dans $S_{\epsilon_n}$ convergeant vers une courbe 
enti\`ere $L$. Celle-ci est dans $S_1\cup S_2$, 
donc par exemple dans $S_1$.
Par hypoth\`ese, $L$ ne peut \^etre 
contenue dans l'intersection 
$S_1 \cap S_2$. Donc, si $L$ coupe $S_2$ hors de $S_3$, par persistance 
d'intersection $L_n$ en fait 
autant. Or $S_{\epsilon_n} \cap S_2$ est contenue dans $S_3$, contradiction.
$\square$

 \null
 
En raffinant cet argument, on s'autorise le contr\^ole d'un seul 
compl\'ementaire :

\thm{Premi\`ere variante} On part des m\^emes trois surfaces dans $\p$ 
avec de plus $d_2=m d_1$. On suppose $S_1\cap S_2$ et
$S_1 \setminus (S_2 \setminus S_3)$ hyperboliques. Alors 
la surface $S_\epsilon=(f_1(f_1^m+\epsilon_1 
f_2)=\epsilon_2 f_3)$ est hyperbolique pour $\epsilon_1, 
\epsilon_2$ petits non 
nuls.

\fthm

En effet, posons $S=
(f_1^m+\epsilon_1 
f_2=0)$. Par ouverture, l'hypoth\`ese
 entra\^ \i ne  l'hyperbolicit\'e de 
$S\setminus (S_2 \setminus S_3)$ pour $\epsilon_1$ 
petit non nul. 
Or $S \cap S_2 = S\cap S_1=S_1 \cap S_2$. 
Ainsi $S \setminus (S_2 \setminus S_3)=
S \setminus (S_1 \setminus S_3)$ et 
$S_1 \setminus (S_2 \setminus S_3)=S_1 \setminus (S 
\setminus S_3)$. Le principe de filtrage s'applique donc 
aux surfaces $S_1$, $S$ et $S_3$. $\square$

\null

En voici un analogue relatif laiss\'e au lecteur (m\^emes notations) :
 
 \thm{Deuxi\`eme variante} Soit de plus $U$ un ouvert d'une 
 quatri\`eme surface $S_4$ en position g\'en\'erale par rapport aux 
 trois premi\`eres. On suppose  $S_1\cap S_2 \setminus U$, 
 $S_1 \cap S_4 
 \setminus S_2$ et $S_1 \setminus ((S_2 \setminus S_3) \cup U)$ 
 hyperboliques. Alors le compl\'ementaire $S_\epsilon \setminus U$ 
 est hyperbolique pour $\epsilon_1,\epsilon_2$ petits non nuls.
 \fthm
 
 La combinaison de ces deux variantes fournit la 
 r\'eduction voulue :
 
 \thm{Lemme} Soient $P$ un plan, $Q,Q'$ deux quadriques et 
 $C,C'$ deux cubiques, en position g\'en\'erale dans $\p$, d'\'equations 
 respectives $(p,q,q',c,c'=0)$. On suppose le compl\'ementaire
 $P \setminus (Q^*\cup C^* )$ hyperbolique, avec 
 $Q^*=Q \setminus C'$ et $C^*=C \setminus Q'$.
 Alors la sextique $S=(c_\epsilon 
(c_\epsilon+\epsilon_3c)=\epsilon_4(q')^3)$, o\`u
 $c_\epsilon=(p(p^2+\epsilon_1q)-\epsilon_2c'$, est hyperbolique 
 pour $\epsilon_1,\epsilon_2,\epsilon_3,
\epsilon_4$ petits non nuls.
\fthm

En effet, notons $C_ \epsilon$ la cubique d'\'equation $(c_\epsilon 
=0)$. 
La deuxi\`eme variante entra\^\i ne que $C_{\epsilon} 
\setminus C^* $ est hyperbolique pour 
$\epsilon_1,\epsilon_2$ petits 
non nuls : en 
prenant $S_1=P, S_2=Q, S_3=C'$ et $U=C^*$, 
il suffit de v\'erifier 
l'hyperbolicit\'e de $P\cap Q \setminus C$ et $P \cap C \setminus Q$;
or $P \cap Q \cap C$ contient $6$ points, et toute courbe 
irr\'eductible priv\'ee de $3$ points est hyperbolique.
Puis, par la premi\`ere variante appliqu\'ee \`a $S_1=C_{\epsilon}, 
S_2=C$ et $S_3=3Q'$, on obtient bien que $S$ est 
hyperbolique puisque $C_{\epsilon} \cap C$ est une 
courbe de genre $10$. $\square$

\null

\thm{Fabrication du filtre}
\fthm
Reste maintenant \`a satisfaire l'hypoth\`ese du lemme. C'est une 
question dans le plan projectif $\pp$.
 Commen\c cons par produire
une conique $Q$ et une cubique $C$ en position g\'en\'erale avec
$\pp \setminus(Q \cup 
C)$ hyperbolique. Partons pour cela d'une configuration 
de deux droites $D_1, D_2$ et d'une cubique $C$,
en 
position g\'en\'erale. Le 
compl\'ementaire $\pp \setminus (D_1 \cup D_2 \cup C)$ est alors
hyperbolique 
d'apr\`es [\BER]. Un argument 
d'ouverture laiss\'e au lecteur 
permet de d\'eformer un peu $D_1 \cup D_2$ en une 
conique $Q$, en pr\'eservant l'hyperbolicit\'e de $\pp \setminus(Q \cup 
C)$.  

 \null
 
Puis on va percer les trous du filtre sur $Q$ et $C$. 
Cela se fait en deux temps. On perturbe d'abord 
de mani\`ere non g\'en\'erique $Q,C$ en $Q'', C''$. On cr\'ee un 
ensemble $T$ de $6$ trous dans $Q\cup C$ d\'efini par $Q\cup C \setminus 
T=(Q\setminus C'') \cup (C \setminus Q'')$, en 
pr\'eservant
l'hyperbolicit\'e de $\pp \setminus (Q \cup C \setminus T)$.  Ensuite on 
perturbe g\'en\'eriquement $Q'',C''$ en $Q',C'$. On obtient maintenant
les $12$ trous voulus, en pr\'eservant l'hyperbolicit\'e de
$\pp \setminus (Q^* \cup C^*)$ (notations du lemme).  
 
 \null

 D\'etaillons la premi\`ere \'etape (la seconde, analogue, sera 
 laiss\'ee au lecteur).
  
 Notons $x_1,\ldots,x_6$ les 6 points de l'intersection $Q \cap C$. On 
 perturbe $Q$ en une conique $Q''$ en conservant les trois premi\`eres 
 intersections avec $C$. On aura donc 
 $$Q'' \cap C=\{x_1,x_2,x_3,t_4,t_5,t_6\}.$$
 Sym\'etriquement, on perturbe $C$ en une cubique $C''$ en conservant 
 les trois derni\`eres intersections avec $Q$, donc $$C'' \cap Q  
 =\{x_4,x_5,x_6,t_1,t_2,t_3\}.$$ Ici $T=\{t_1,\ldots,t_6\}$ 
 est disjoint de $Q \cap C$. Montrons que, si $Q'',C''$ sont assez 
 proches de $Q,C$, le compl\'ementaire  
 $\pp \setminus (Q \cup C 
\setminus T)$ est encore hyperbolique.
 
  \null

Sinon, on produit une suite $(L_n)$ de courbes enti\`eres dans 
$\pp \setminus
(Q \cup  C\setminus T_n)$ (avec des notations \'evidentes) convergeant
vers une 
courbe enti\`ere $L$. Remarquons que la courbe $L_n$ \'evite $C$ 
pr\`es de $x_1,x_2,x_3$, et \'evite $Q$ pr\`es de 
$x_4,x_5,x_6$. 

V\'erifions que $L$ n'est pas contenue dans $Q \cup C$ : en 
effet, si 
elle \'etait contenue dans $C$ par exemple, elle couperait $Q$ en l'un des 
points $x_4,x_5,x_6$ ($C$ priv\'ee de trois points 
est hyperbolique); par persistance d'intersection $L_n$ 
couperait 
encore $Q$ 
pr\`es de l'un de ces points, contradiction.

Maintenant $L$ doit couper $Q \cup C$ puisque 
$\pp \setminus (Q \cup C)$ est hyperbolique. Par persistance 
d'intersection et comme $T_n$ converge vers 
$Q\cap C$, $L$ coupe en fait $Q \cup 
C$ en un point de $Q\cap C$. Donc, toujours par persistance d'intersection, 
$L_n$ va 
couper \`a la fois $C$ et $Q$ pr\`es d'un point de $Q\cap C$, 
contradiction. $\square$

 \Refs

\widestnumber\no{99}
\refno=0
  
\bref \by F. Berteloot, J. Duval \paper Sur l'hyperbolicit\'e de 
certains compl\'ementaires \jour Ens. Math. 
\vol47\yr2001\pages253--267
\endref
 
\bref \by R. Brody \paper Compact manifolds and hyperbolicity \jour 
Trans. Amer. Math. Soc. \vol235\yr1978\pages213--219
\endref

\bref \by J.-P. Demailly, J. El Goul \paper Hyperbolicity of generic 
surfaces of high degree in projective 3-space \jour Amer. J. Math. 
\vol122\yr2000\pages515--546
  \endref
  
  \bref \by B. Shiffman, M. Zaidenberg \paper New examples of 
  hyperbolic octic surfaces in $P^3$ \jour preprint 2003 
  arXiv math.AG/0306360
   \endref
  
  \bref \by M. Zaidenberg \paper Hyperbolic surfaces in $P^3$: 
  examples \jour preprint 2003 arXiv math.AG/0311394
  \endref
  
\endRefs

\address 
\noindent  
Laboratoire \'Emile Picard, UMR CNRS 5580, 
  Universit\'e Paul Sabatier, 31062 Toulouse Cedex 4, France.
 \endaddress
\email 
  duval\@picard.ups-tlse.fr 
\endemail

\enddocument